\documentclass[12pt]{article}

\usepackage{amsmath}
\usepackage{amssymb}
\usepackage{amsfonts}
\usepackage{amsmath}
\usepackage{amsthm}
\usepackage{tikz}
\usepackage{pgf}
\usepackage[bookmarks=true,pdfborder={0 0 0}]{hyperref}

\newtheorem{theorem}{Theorem}[section]
\newtheorem{lemma}[theorem]{Lemma}
\newtheorem{corollary}[theorem]{Corollary}
\newtheorem{prop}[theorem]{Proposition}

\theoremstyle{definition}

\theoremstyle{remark}

\numberwithin{equation}{section}

\title{Alternative Proofs on the Indices of Cacti and Unicyclic Graphs with $n$ Vertices}
\author{Sudipta Mallik\\ \normalsize{Department of Mathematics}\\
\normalsize{University of Wyoming, Laramie, WY, USA}\\
\normalsize{E-mail: \href{mailto:smallik@uwyo.edu}{smallik@uwyo.edu}}}

\begin{document}
\maketitle

\begin{center}
\vspace{10pt}
\end{center}



\begin{abstract}
 Let $H_n$ be the cactus obtained from the star $K_{1,n-1}$ by adding $\left \lfloor \frac{n-1}{2} \right \rfloor$ independent edges between pairs of pendant vertices.  Let $K_{1,n-1}^+$ be the unicyclic graph obtained from the star $K_{1,n-1}$ by appending one edge. In this paper we give alternative proofs of the following results: Among all cacti with $n$ vertices,  $H_n$ is the unique cactus whose spectral radius is maximal, and among all unicyclic graphs with $n$ vertices, $K_{1,n-1}^+$  is the unique unicyclic graph whose spectral radius is maximal. We also prove that among all odd-cycle graphs with $n$ vertices,  $H_n$ is the unique odd-cycle graph whose spectral radius is maximal. \end{abstract}

\section{Introduction}
Let $G$ be a simple graph with vertex set $\{v_1, v_2, . . . , v_n\}$. The adjacency matrix of $G$, $A(G) = [a_{ij}]$ is defined to be the $n \times n$ matrix such that  $a_{ij}= 1$ if $v_i$ is adjacent to $v_j$, and $ a_{ij}= 0$ otherwise. Since $A(G)$ is symmetric, all of its eigenvalues are real. The spectral radius of $G$, $\rho (G)$, is the largest eigenvalue of  $A(G)$ and  it is also called the {\it index} of  $G$.  

When $G$  is connected, $A(G)$ is irreducible and by the Perron-Frobenius Theorem [5, p 181], $\rho (G)$ is a simple eigenvalue of $A(G)$ and there is a unique positive unit eigenvector corresponding to $\rho (G)$. This eigenvector is called the $\it{Perron \;vector}$ of $G$.

A {\it pendant} vertex is a vertex of degree 1. We call an edge a {\it pendant edge} if it is a bridge connecting a pendant vertex. Lets denote the degree of a vertex $v$ by $d(v)$. Let $\Delta (G)$ denote the highest degree of all vertices of $G$. We denote   the set of all vertices adjacent to $v$ by $N(v)$. 


A path is called an {\it odd-path} if its length (i.e., the number of its edges) is odd. Otherwise it is called an {\it even-path}. A cycle is called an {\it odd-cycle} if its length  (i.e., the number of its edges) is odd. Otherwise it is called an {\it even-cycle}. A graph is called an {\it odd-cycle graph} if each of its cycles is an {\it odd-cycle}.

A graph is called a {\it cactus} if its cycles have at most one common vertex. Let $H_n$ be the cactus obtained from the star $K_{1,n-1}$ by adding $\left \lfloor \frac{n-1}{2} \right \rfloor$ independent edges between pairs of pendant vertices (see Fig. 2).
Let $\mathcal C(n)$ be the set of all cacti of order $n$ (i.e., with $n$ vertices).

A connected graph with a unique cycle is called a {\it unicyclic graph}. So a unicyclic graph can be seen as a tree with an extra edge. By $K_{1,n-1}^+$ we denote the unicyclic graph   obtained from the star $K_{1,n-1}$ by appending one edge (see Fig. 3). Let $\mathcal U(n)$ be the set of all unicyclic graphs of order $n$. 

\section{Main Results}
\begin{theorem} $ [2, \mbox{Thm.}1]$  Let $u,v$ be two vertices of the connected graph $G$. Suppose $v_1, v_2,\ldots,v_s\;(1\leq s \leq d(v))$ are some vertices of $N(v)\backslash N(u)$ and $x=(x_1, x_2,\ldots,x_n)^T$ is the Perron vector of $G$, where $x_i$ corresponds to the vertex $v_i\;(1 \leq i \leq n)$. Let $G^*$ be the graph obtained from $G$ by deleting the edges $vv_i$ and adding the edges $uv_i \;(1 \leq i \leq s)$. If $x_u\geq x_v$, then $\rho (G^*)>\rho (G)$.\\ 
\end{theorem}

\begin{center}
\begin{tikzpicture}
[colorstyle/.style={circle, draw=black!100,fill=black!100, thick, inner sep=0pt, minimum size=.5 mm},>=stealth]
\node at (0,1)[colorstyle, label=right:$v_1$]{};
\node at (0,.75)[colorstyle, label=right:$v_2$]{};
\node at (0,.5)[colorstyle]{};
\node at (0,.25)[colorstyle]{};
\node at (0,0)[colorstyle, label=right:$v_s$]{};
\node at (-2,.75)[colorstyle, label=left:$v$]{};
\node  at (-2,0)[colorstyle, label=left:$u$]{};
\draw [thick] (0,1)--(-2,.75)--(0,0);
\draw [thick] (-2,.75)--(0,.75);

\node at (-1,-.5)[label=below:$G$]{};

\node at (2,0)[colorstyle, label=left:$u$]{};
\node  at (2,.75)[colorstyle, label=left:$v$]{};
\node  at (4,0)[colorstyle, label=right:$v_s$]{};
\node at (4,.25)[colorstyle]{};
\node at (4,.5)[colorstyle]{};
\node at (4,.75)[colorstyle, label=right:$v_2$]{};
\node at (4,1)[colorstyle, label=right:$v_1$]{};
\draw [thick] (4,1)--(2,0)--(4,.75);
\draw [thick] (2,0)--(4,0);
\node at (3,-.5)[label=below:$G^*$]{};

\end{tikzpicture}\\
\it Figure 1
\end{center}

Using mainly this theorem  we will give an alternative proof of the following theorem:
\begin{theorem} $[3, \mbox{Thm. } 3.1]$ Let $G \in \mathcal C(n)$. Then $\rho (G)\leq \rho (H_n)$, equality holds if and only if $G \cong H_n$. \end{theorem}

\begin{center}
\begin{tikzpicture}
[colorstyle/.style={circle, draw=black!100,fill=black!100, thick, inner sep=0pt, minimum size=.5 mm},>=stealth]
\node at (-.75,1)[colorstyle, label=above:$v_2$]{};
\node at (-1.5,1)[colorstyle, label=above:$v_3$]{};
\node at (-2.25,1)[colorstyle, label=above:$v_4$]{};
\node at (-3,1)[colorstyle, label=above:$v_5$]{};
\node at (-3.25,1)[colorstyle]{};
\node at (-3.5,1)[colorstyle]{};
\node at (-3.75,1)[colorstyle]{};
\node at (-4.25,1)[colorstyle, label=above:$v_{n-1}$]{};
\node at (-5.,1)[colorstyle, label=above:$v_{n}$]{};

\node  at (-3,-1)[colorstyle, label=below:$v_1$]{};
\draw [thick] (-3,-1)-- (-.75,1)--(-1.5,1)--(-3,-1)--(-2.25,1)--(-3,1)--(-3,-1)--(-4.25,1)--(-5,1)--(-3,-1);
\node at (-3,-1.5)[label=below:$n$ is odd]{};

\node at (5,1)[colorstyle, label=above:$v_2$]{};
\node at (4.25,1)[colorstyle, label=above:$v_3$]{};
\node at (3.5,1)[colorstyle, label=above:$v_4$]{};
\node at (2.75,1)[colorstyle, label=above:$v_5$]{};
\node at (2.5,1)[colorstyle]{};
\node at (2.25,1)[colorstyle]{};
\node at (2,1)[colorstyle]{};
\node at (1.5,1)[colorstyle, label=above:$v_{n-2}$]{};
\node at (.75,1)[colorstyle, label=above:$v_{n-1}$]{};
\node at (.75,-1)[colorstyle, label=below:$v_{n}$]{};
\node  at (2.25,-1)[colorstyle, label=below:$v_1$]{};
\draw [thick] (2.25,-1)-- (5,1)-- (4.25,1)--(2.25,-1)--(3.5,1)--(2.75,1)-- (2.25,-1)--(1.5,1)--(.75,1)-- (2.25,-1)--(.75,-1);
\node at (-3,-1.5)[label=below:$n$ is odd]{};

\node at (2.5,-1.5)[label=below:$n$ is even]{};

\end{tikzpicture}\\
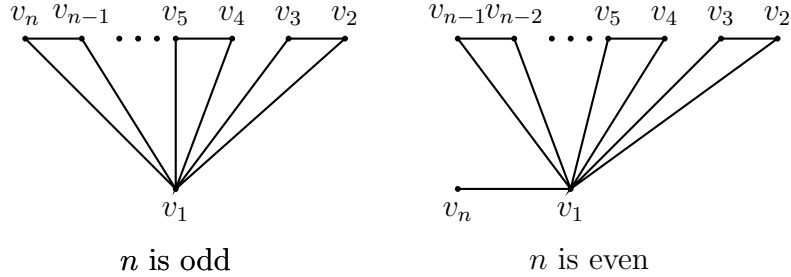
\it Figure 2: $H_n$
\end{center}

First we prove the above theorem for all  connected cacti with maximal number of edges as following.
\begin{theorem} Let $G \in \mathcal C(n)$ be connected with maximal number of edges. Then $\rho (G)\leq \rho (H_n)$, equality holds if and only if $G \cong H_n$. \end{theorem}

Before proving the theorem we will  record the following  propositions regarding connected cactus $G$ of order $n$ with maximal number of edges:
\begin{prop} Let $G$ be a connected cactus $G$ of order $n$ with maximal number of edges. Then\\ (a) All cycles of $G$ are triangles with at most one edge not in some triangle except when $G \cong C_4$.[1,  Lemma 6.7]\\
(b) If $n \leq 5$ and $G \ncong C_4$ then $G \cong H_n$.\\
(c) Let $t(G):= \mbox{the number of vertices of $G$ of degree} \geq 3$. If $n \geq 6$ then $t(G)=1$ if and only if $G \cong H_n$.\\
(d) If $t(G)>1$ then there are always two adjacent vertices of degree $\geq 3$.[from (a)]\end{prop}

\lemma Let $G$ be a connected cactus of order $n$ with maximal number of edges. Let $u$ and $v$ be two adjacent vertices of $G$ of degree $\geq 3$ such that $N(v)\backslash \{u \cup N(u)\} = \{v_1, v_2,\ldots,v_{s}\}$.  Let $G_1$ be the graph obtained from $G$ by deleting the edges $vv_i$ and adding the edges $uv_i \;(1 \leq i \leq s)$. Then $G_1$ is also a connected cactus of order $n$ with maximal number of edges.
\proof Since $u$ and $v$ are adjacent  and $G$ is connected then $G_1$ is connected. $G_1$ and $G$ have the same number of edges since in making of $G_1$ the numbers of deleted edges and added edges are same. So it suffices to show that $G_1$ is a cactus. By Proposition 2.4(a), all cycles of $G$ are triangles with at most one edge not in some triangle. If $s$ is even $v_1, v_2,\ldots,v_{s}$ form
exactly $\frac{s}{2}$ triangles at $v$ having no other common vertex. These $\frac{s}{2}$ triangles at $v$ corresponds $\frac{s}{2}$ branches of  induced subgraphs of $G$ which are also connected cacti with maximal number of edges having a unique common vertex $v$. When we delete $vv_i$ and add  $uv_i \;(1 \leq i \leq s)$, $v$ becomes a vertex of degree $\leq 2$ and $u$ is added with  $\frac{s}{2}$ branches of  connected cacti with maximal number of edges having a unique common vertex $u$. So $G_1$ is a  connected cactus with maximal number of edges. If $s$ is odd $v_1, v_2,\ldots,v_{s}$ form  one edge and exactly $\frac{s-1}{2}$ triangles at $v$ having no other common vertex. Then  by similar arguments  $G_1$ is a connected cactus with maximal number of edges.\qed\\

{ \it Proof of Theorem 2.3.}  Suppose $G$ is a connected cactus of order $n$ with maximal number of edges such that $\rho (G) \geq \rho (G')$ for all connected cactus $G'$ of order $n$ with maximal number of edges.
If $G \cong H_n$ then there is nothing to prove. Let $G \ncong H_n$. Now $G \ncong C_4$ since $\rho (H_4) > \rho (C_4)$. Then by Proposition 2.4(c), $t(G)>1$. Now by Proposition 2.4(d),  suppose $u$ and $v$ are two adjacent vertices of $G$ of degree $\geq 3$.

Let $x=(x_1, x_2,\ldots,x_n)^T$ be the Perron vector of $G$, where $x_i$ corresponds to the vertex $v_i\;(1 \leq i \leq n)$. Suppose  $x_u\geq x_v$. 
By Proposition 2.4(a), $G$ has no 4-cycle. So $u$ and $v$ can have at most one common adjacent vertex. Since $v$  has degree at least 3, $N(v)\backslash \{u \cup N(u)\} \neq \phi$. Let $v_1, v_2,\ldots,v_{s}$ be all vertices  of $N(v)\backslash \{u \cup N(u)\}$.  Let $G_1$ be the graph obtained from $G$ by deleting the edges $vv_i$ and adding the edges $uv_i \;(1 \leq i \leq s)$. By Lemma 2.5,  $G_1$ is also a connected cactus of order $n$ with maximal number of edges. Now by Theorem 2.1, $\rho (G_1)>\rho (G)$ which is a contradiction to the fact that $\rho (G) \geq \rho (G')$ for all edge maximal connected cactus $G'$ of order $n$.\qed \\

\corollary  Let $G$ be a connected cactus of order $n$ with maximal number of edges such that $t(G)>1$. Then there exists a connected cactus $G_1$ of order $n$  with maximal number of edges, not necessarily isomorphic to $H_n$, such that $t(G_1)=t(G)-1$ and $\rho (G_1)>\rho (G)$.
\proof From the proof of Lemma 2.5 it is clear that  when we delete  $vv_i$ and add $uv_i \;(1 \leq i \leq s)$ in $G$, $v$ becomes a vertex  in $G_1$ of degree $\leq 2$. So  $t(G_1)=t(G)-1$. Now by the above proof $\rho (G_1)>\rho (G)$.\qed\\

Let $G -kl$ denote the graph $G$ without the edge $kl$. 
\begin{lemma}$[1, \mbox{Lemma\;} 6.4]$ $\rho (G) \geq \rho (G-kl)$ for any edge $kl$ of $G$, with strict inequality when $G$ is connected.\end{lemma}
 By the above lemma  it suffices to prove Theorem 2.2 for edge maximal connected cacti. Since we already proved Theorem 2.2 for connected cacti with maximal number of edges in Theorem 2.3, then to prove Theorem 2.2 it suffices to prove the following theorem.
\theorem Let $G \in \mathcal C(n)$ be edge maximal connected. Then there is a connected cactus $G^*$ with maximal number of edges such that $\rho (G)\leq \rho (G^*)$.
\proof If $G$ is a connected cactus with maximal number of edges there is nothing to prove. Suppose  $G$ is an edge maximal connected cactus without maximal number of edges. Let $C$ be a cycle of $G$ and $e$ be an edge at a vertex of $C$ but not in $C$. Since $G$ is edge maximal then one of the following is true.\\
(a) $e$ is in a cycle $C' \neq C$ in $G$.\\
(b) $e$ is a bridge between two cycles $C$ and $C' \neq C$ in $G$.\\
(c) $e$ is a pendant edge in $G$.\\
Since $G$ is edge maximal there are no two consecutive edges which are not in any circle. Because if $uv$ and $vw$ are such two, then we can add a new edge $uw$ while the new graph is still a connected cactus. Now  we construct $G^*$ from $G$ using the following steps.\\
Step 1. Now let $C_k$, $k\geq 4$ be a cycle in $G$. Let $u$ and $v$ be two adjacent vertices in $C_k$. Since $k\geq 4$, suppose $u$ is adjacent to $x (\neq v)$ and $v$ is adjacent to $y (\neq u)$. Then $N(u)\backslash \{v \cup N(v)\}=\{x\}$ and  $N(v)\backslash \{u \cup N(u)\}=\{y\}$ in $C_k$. Let $x=(x_1, x_2,...,x_n)^T$ be the Perron vector of $G$, where $x_i$ corresponds to the vertex $v_i\;(1 \leq i \leq n)$. Suppose  $x_u\geq x_v$. Then deleting the edge $vy$ and adding the edge $uy$ in $G$ we will get a graph $G'$ which is same as $G$ except in which $C_k$ becomes $C_{k-1}$ joined with the edge $uv$. Now by Theorem 2.1, $\rho (G')> \rho (G)$. Repeating this process in every $C_k$, $k\geq 4$ in $G$ we get  a connected cactus $G_1$ of order $n$  in which cycles are triangles and $\rho (G_1)> \rho (G)$. If $G_1$ has at most one edge not in any triangle then by  Proposition 2.4(a), $G_1$ is a connected cactus with maximal number of edges. Then we are done.

Step 2. Suppose $G_1$ has at least two edges not in any triangle. If $uv$ and $vw$ are such two, then we can add a new edge $uw$ producing an extra triangle. Repeating this in all possible cases we get  a connected cactus $G_2$ of order $n$  in which cycles are triangles and by Lemma 2.7, $\rho (G_2)\geq \rho (G_1)$, equality holds if and only if there are no two consecutive edges that are not in any cycle. Similarly if $G_2$ has at most one edge not in any triangle then by  Proposition 2.4(a), $G_1$ is a connected cactus with maximal number of edges. Then we are done.

Step 3.  Suppose $G_2$ has at least two edges not in any triangle. By construction of $G_2$ there are no  two consecutive edges which are not in any triangle. Let $uv$ be an edge that is not in any triangle such that $uv$ is a bridge between two triangles in $G_2$. Then obviously $N(u)\backslash \{v \cup N(v)\}\neq \phi$ and  $N(v)\backslash \{u \cup N(u)\}\neq \phi$.  Let $x=(x_1, x_2,\ldots,x_n)^T$ be the Perron vector of $G_2$, where $x_i$ corresponds to the vertex $v_i\;(1 \leq i \leq n)$. Suppose $x_u\geq x_v$.  Let $N(v)\backslash \{u \cup N(u)\}=\{v_1,v_2,\cdots, v_s\}$.  Let $G_2^*$ be the graph obtained from $G_2$ by deleting the edges $vv_i$ and adding the edges $uv_i \;(1 \leq i \leq s)$. Then $G_2^*$ is a connected cactus of order $n$  in which cycles are triangles with the pendant edge $uv$ and by Theorem 2.1, $\rho (G_2^*)> \rho (G_2)$.  Repeating this process for all bridges between triangles we can get a connected cactus $G_3$ of order $n$   in which cycles are triangles and by Theorem 2.1,  $\rho (G_3)> \rho (G_2)$. If $G_3$ has at most one edge not in any triangle then by  Proposition 2.4(a), $G_3$ is a connected cactus with maximal number of edges. Then we are done. 

Step 4.   Suppose $G_3$ has at least two edges not in any triangle. By construction of $G_3$, all the edges  of $G_3$ which are not in any triangle are pendant edges. Let $ux$ and $vy$ be two pendant edges of $G_3$. If $x=y$, add the edge $uv$ to $G_3$ which increases number of triangle and also spectral radius of $G_3$ by Lemma 2.7. We will do this in all possible cases and get a connected cactus $G_3^*$ in which cycles are triangles and by Lemma 2.7, $\rho (G_3^*)> \rho (G_3)$. Now let $x\neq y$. Then $N(u)\backslash \{v \cup N(v)\}=\{x\}$ and  $N(v)\backslash \{u \cup N(u)\}=\{y\}$.  Let $x=(x_1, x_2,\ldots,x_n)^T$ be the Perron vector of $G_3^*$, where $x_i$ corresponds to the vertex $v_i\;(1 \leq i \leq n)$. Suppose  $x_u\geq x_v$. Then deleting the edge $vy$ and adding the edge $uy$ in $G$ we will get a graph $G_4$ of order $n$. Then $G_4$ is a connected cactus of order $n$  in which cycles are triangles and by Theorem 2.1, $\rho (G_4)> \rho (G_3)$. Repeating this process for all pendant edges we can get a connected cactus $G_5$ of order $n$  in which cycles are triangles and by Theorem 2.1,  $\rho (G_5)> \rho (G_4)$.

Step 5. By construction of $G_5$, all the pendant edges  of $G_5$,  except at most one, form pairs having a common vertex. Now joining corresponding two pendant vertices of each pair we can form a new triangle for each such pair. Then we  get a connected cactus $G^*$ of order $n$   in which cycles are triangles and by Lemma 2.7,  $\rho (G^*)> \rho (G_5)$. By construction of $G^*$, it has at most one pendant edge. Then $\rho (G^*)\geq \rho (G)$ and by  Proposition 2.4(a), $G^*$ is a connected cactus of order $n$ with maximal number of edges.\qed \\

Now lets prove a corollary of Theorem 2.2 as following.

\begin{corollary}  For all odd-cycle graph $G$, $\rho (G)\leq \rho (H_n)$, equality holds if and only if $G \cong H_n$. \end{corollary}

Before proving this corollary we will prove the following lemma.

\lemma Every odd-cycle graph is a cactus.
\proof  Let $G$ be an odd-cycle graph. Suppose $G$ is not a cactus. Then  $G$  have two odd cycles, say $C$ and $C'$ such that they have at least two common vertices. Let  $v_1, v_2,\ldots,v_{k}$ be all common vertices of $C$ and $C'$. So these vertices divide each of $C$ and $C'$ into a series of consecutive paths, say $P^{1}, P^{2},\ldots, P^{k}$ for $C$ and $P^{1'}, P^{2'},\ldots, P^{k'}$ for $C'$ where $P^{i}$ is the path from $v_i$ to $v_{i+1}$ in $C$ and  $P^{i'}$ is the path from $v_i$ to $v_{i+1}$ in $C'\;(1 \leq i \leq k)$, assuming $v_{k+1}=v_1$ . Since $C\neq C'$, then $P^i \neq P^{i'}$ for some $i$. If $P^i$ and $P^{i'}$ both are even-paths or odd-paths then $P^i \cup P^{i'}$ is an even-cycle in $G$ - a contradiction. Otherwise suppose $P^i$ is an even-path and $P^{i'}$ is an odd-path. Let $P$ be a path from $v_i$ to $v_{i+1}$ obtained from $C$ by deleting nonpendant vertices of $P^i$ and corresponding incident edges. Since $C$ is an odd-cycle and $P^i$ is an even-path, then $P$ is an odd-path from  $v_i$ to $v_{i+1}$  in $G$. Now since odd-paths $P$ and $P^{i'}$ are disjoint except at the end points $v_i$ and $v_{i+1}$, then $P\cup P^{i'}$ is an even-cycle in $G$ - a contradiction. \qed\\

{\it Proof of Corollary 2.9.}  From the above lemma every odd-cycle graph is a cactus. So a connected odd-cycle graph is a connected cactus. Now by Theorem 2.2 for all connected cactus $G$,  $\rho (G)\leq \rho (H_n)$, equality holds if and only if $G \cong H_n$. Since $H_n$ is a connected odd-cycle graph the corollary follows.\qed\\

\begin{theorem} $[4]$  Let $G \in \mathcal U(n)$. Then $\rho (G)\leq \rho (K_{1,n-1}^+)$, equality holds if and only if $G \cong K_{1,n-1}^+$. \end{theorem}

\begin{center}
\begin{tikzpicture}
[colorstyle/.style={circle, draw=black!100,fill=black!100, thick, inner sep=0pt, minimum size=.5 mm},>=stealth]
\node at (-.75,1)[colorstyle, label=above:$v_2$]{};
\node at (-1.5,1)[colorstyle, label=above:$v_3$]{};
\node at (-2.25,1)[colorstyle, label=above:$v_4$]{};
\node at (-3,1)[colorstyle, label=above:$v_5$]{};
\node at (-3.25,1)[colorstyle]{};
\node at (-3.5,1)[colorstyle]{};
\node at (-3.75,1)[colorstyle]{};
\node at (-4.25,1)[colorstyle, label=above:$v_{n-1}$]{};
\node at (-5.,1)[colorstyle, label=above:$v_{n}$]{};

\node  at (-3,-1)[colorstyle, label=below:$v_1$]{};
\draw [thick] (-3,-1)-- (-.75,1)--(-1.5,1)--(-3,-1)--(-2.25,1);
\draw [thick](-3,1)--(-3,-1)--(-4.25,1);
\draw [thick](-5,1)--(-3,-1);

\end{tikzpicture}\\
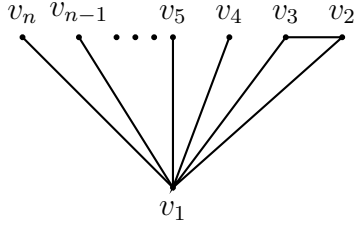
\it Figure 3: $K_{1,n-1}^+$
\end{center}

Using mainly Theorem 2.1  we will give an alternative proof of the above theorem. Before proving the theorem we will record the following  propositions regarding unicyclic graph  $G$ of order $n$ :
\begin{prop}  (a) Let $G \in \mathcal U(n)$.  Two adjacent vertices in $G$ have at most one common vertex. Two vertices in $G$ have one common vertex if and only if they are in a triangle.\\
(b)  Let $G \in \mathcal U(n)$. If $n=3$ then  $G \cong K_{1,n-1}^+$.\\
(c)  Let $G \in \mathcal U(n)$. $G \cong K_{1,n-1}^+$ if and only if $\Delta (G)=n-1$.\\
(d) Let $G$ be a graph of order $n$. Then $G$ is  unicyclic  if and only if $G$ is connected having exactly $n$ edges.\\
\end{prop}




\lemma Let $G \in \mathcal U(n)$ and  $G \ncong K_{1,n-1}^+$. Then there are two adjacent vertices $u$ and $v$ in $G$ such that $N(u)\backslash \{v \cup N(v)\} \neq \phi$ and  $N(v)\backslash \{u \cup N(u)\}\neq \phi$. 
\proof  
Since $G \ncong K_{1,n-1}^+$  then by Proposition 2.12(b), $n\geq 4$. If $G\cong C_n$, $n \geq 4$ then $d(v)=2$ for every vertex $v$ of $G$. Let $u$ and $v$ be two adjacent vertices in $G \cong C_n$, $n \geq 4$. Suppose $u$ is adjacent to $x (\neq v)$ and $v$ is adjacent to $y (\neq u)$. Then $N(u)\backslash \{v \cup N(v)\}=\{x\}$ and  $N(v)\backslash \{u \cup N(u)\}=\{y\}$. Suppose $G\ncong C_n$, $n \geq 4$. Then there is a vertex $v$ in $G$ such that $d(v)\geq 3$.  Since $G \ncong K_{1,n-1}^+$  by Proposition 2.12(c),  $d(v)\leq \Delta (G)<n-1$. So $v$ is not adjacent to at least one vertex in $G$. Let $w$ be one such. Suppose $P$ is a shortest path between $v$ and $w$.
Take the vertex  adjacent to $v$ in $P$ as $u$.   Now by Proposition 2.12(a), $v$ and $u$ have at most one common vertex, say $x$. Since $d(v) \geq 3$, $v$ is adjacent to a vertex, say $y$ that is different from  $u$ and $x$. Then $y \in N(v)\backslash \{u \cup N(u)\}$. Let  $z$ be the other adjacent vertex of $u$ in $P$. If $z$ is not adjacent to $v$ then $z \in N(u)\backslash \{v \cup N(v)\}$. Otherwise we have a path from $v$ to $w$ of length shorter than that of $P$, a contradiction.\qed\\

\lemma  Let $G \in \mathcal U(n)$ and  $G \ncong K_{1,n-1}^+$. Let $u$ and $v$ be two adjacent vertices in $G$ such that $N(u)\backslash \{v \cup N(v)\} \neq \phi$ and $N(v)\backslash \{u \cup N(u)\}\neq \phi$. 

1. Let $N(u)\backslash \{v \cup N(v)\}=\{u_1,u_2,\ldots,u_t\}$. Let $G_1$ be the graph obtained from $G$ by deleting the edges $uu_i$ and adding the edges $vu_i \;(1 \leq i \leq t)$. Then $G_1$ is also a unicyclic graph of order $n$.

2. Let $N(v)\backslash \{u \cup N(u)\}=\{v_1, v_2,\ldots,v_{s}\}$.   Let $G_1$ be the graph obtained from $G$ by deleting the edges $vv_i$ and adding the edges $uv_i \;(1 \leq i \leq s)$. Then $G_1$ is also a unicyclic graph of order $n$.
\proof 1. Since $u$ and $v$ are adjacent and $G$ is connected, then $G_1$ is also connected. Also it is clear that $G$ and $G_1$ have same number of edges which is $n$. So $G_1$ is a connected graph of order $n$ having exactly $n$ edges. Then by Proposition 2.12(d), $G_1$ is unicyclic.

2.  It follows from similar arguments.\qed\\

{ \it Proof of Theorem 2.11.}  Let $G \in \mathcal U(n)$ such that $\rho (G) \geq \rho (G')$ for all $G' \in \mathcal U(n)$. If $G \cong K_{1,n-1}^+$ there is nothing to prove. Let  $G \ncong K_{1,n-1}^+$.  By Lemma 2.13, there are two adjacent vertices $u$ and $v$ in $G$ such that $N(u)\backslash \{v \cup N(v)\} \neq \phi$ and $N(v)\backslash \{u \cup N(u)\} \neq \phi$.   Let $x=(x_1, x_2,\ldots,x_n)^T$ be the Perron vector of $G$, where $x_i$ corresponds to the vertex $v_i\;(1 \leq i \leq n)$. 

Case 1. $x_u\geq x_v$. Let $v_1, v_2,\ldots,v_{s}$ be all vertices  of $N(v)\backslash \{u \cup N(u)\}$.  Let $G_1$ be the graph obtained from $G$ by deleting the edges $vv_i$ and adding the edges $uv_i \;(1 \leq i \leq s)$. By Lemma 2.14,  $G_1$ is also a unicyclic graph of order $n$. Now by Theorem 2.1, $\rho (G_1)>\rho (G)$.

Case 2.  $x_v\geq x_u$.  Let $u_1, u_2,\ldots,u_{t}$ be all vertices  of $N(u)\backslash \{v \cup N(v)\}$.  Let $G_1$ be the graph obtained from $G$ by deleting the edges $uu_i$ and adding the edges $vu_i \;(1 \leq i \leq t)$. By Lemma 2.14,  $G_1$ is also a unicyclic graph of order $n$. Now by Theorem 2.1, $\rho (G_1)>\rho (G)$.
 
In either case $\rho (G_1)>\rho (G)$ which is a contradiction to the fact that  $\rho (G) \geq \rho (G')$ for all $G' \in \mathcal U(n)$.\qed\vspace{20pt}\\
\textbf{Acknowledgements}\vspace{5 pt}\\
The author would like to thank his academic  advisor Bryan Shader for his valuable suggestions.\\

\end{document}